# Modeling On-Line Art Auction Dynamics Using Functional Data Analysis

**Srinivas K. Reddy and Mayukh Dass**

*Abstract.* In this paper, we examine the price dynamics of on-line art auctions of modern Indian art using functional data analysis. The purpose here is not just to understand what determines the final prices of art objects, but also the price movement during the entire auction. We identify several factors, such as artist characteristics (established or emerging artist; prior sales history), art characteristics (size; painting medium—canvas or paper), competition characteristics (current number of bidders; current number of bids) and auction design characteristics (opening bid; position of the lot in the auction), that explain the dynamics of price movement in an on-line art auction. We find that the effects on price vary over the duration of the auction, with some of these effects being stronger at the beginning of the auction (such as the opening bid and historical prices realized). In some cases, the rate of change in prices (velocity) increases at the end of the auction (for canvas paintings and paintings by established artists). Our analysis suggests that the opening bid is positively related to on-line auction price levels of art at the beginning of the auction, but its effect declines toward the end of the auction. The order in which the lots appear in an art auction is negatively related to the current price level, with this relationship decreasing toward the end of the auction. This implies that lots that appear earlier have higher current prices during the early part of the auction, but that effect diminishes by the end of the auction. Established artists show a positive relationship with the price level at the beginning of the auction. Reputation or popularity of the artists and their investment potential as assessed by previous history of sales are positively related to the price levels at the beginning of the auction. The medium (canvas or paper) of the painting does not show any relationship with art auction price levels, but the size of the painting is negatively related to the current price during the early part of the auction. Important implications for auction design are drawn from the analysis.

*Srinivas K. Reddy is Robert O. Arnold Professor of Business and Director, Coca-Cola Center for Marketing Studies, Terry College of Business, University of Georgia, Athens, Georgia 30602, USA e-mail:* sreddy@terry.uga.edu. *Mayukh Dass is a Doctoral Candidate, Terry College of Business, University of Georgia, Athens, Georgia 30602, USA e-mail:* dass@uga.edu.

## 1. INTRODUCTION

The popularity and success of on-line auctions and the availability of detailed, unique bidding data have made it possible for researchers to address interesting and important questions regarding seller and bidder behavior and auction price dynamics, hitherto not possible with traditional auctions. Pioneering work using on-line auction data was done by Bapna and his associates [4, 5, 6, 7] to address important design issues of on-line auctions. Jank, Shmueli and their associates, in a series of papers, introduced innovative statistical methods to visualize and analyze the dynamics of on-line auctions that would not have been possible using traditional auction data and standard statistical tools [8, 15, 26, 27]. In the business to business (B2B) on-line auction space, Elmaghraby [14] has looked at the pricing practices and the challenges faced in B2B e-marketplaces. Jap [16, 17] and Mithas, Jones and Mitchell [20] have shed some light on the buyer–seller behavior in B2B on-line reverse auctions. More recently, Spann and Tellis [29] examined the rationality of consumer decision-making in name-your-own-price auctions.

The success of on-line auctions has sparked innovation in traditional art auction houses such as Christie's and Sotheby's to adopt the new technologies to facilitate easier viewing and bidding. More importantly, a new set of auction houses, particularly in some of the emerging art markets, has surfaced, houses which have fully incorporated/embraced the on-line auction format with great success.

Fine art is considered to be a hedonic product, the consumption of which is driven more by affective experience such as the aesthetic pleasure that one derives from it rather than its utilitarian or functional benefits. Research on art has mostly concentrated on price formation at traditional auctions and art's potential as an investment [1, 2, 9, 19, 23]. (Ashenfelter and Graddy [2] report from their analysis of 15 studies an estimated average real return of art as investment of 2.6% with a range of 0.55% to 5%.) Typically the data have consisted of realized prices at major auction houses over a number of years. Modeling of these data was typically accomplished by standard statistical and econometric methods. The aggregated nature of the data available and the lack of bidding information (such as the opening bid, the bid arrival data, number of bidders, total number of bids, etc.) for specific auctions precluded one from looking at the dynamics of price formation during the auction.

Analysis of on-line auctions of hedonic heterogeneous products like fine art is both complex and challenging. First, the uniqueness of the artwork and the artists provides a modeling challenge due to the high variability across the auctioned lots. This variability is a result of the uniqueness and scarcity of the art objects and the genre differences among the artists. Second, these artworks possess a more subjective private (hedonic) value to the bidders than an objective common value. Finally, a small group of dedicated art collectors and investors follow on-line art auctions. The emergence of on-line art auctions and the availability of bidding dynamic data provide us with a unique opportunity to develop and explore models to analyze price dynamics of such hedonic products.

Unique data offered in on-line auctions require unique methods of analysis. In this paper, we present the application of functional data analysis [25], which is well suited to explore the dynamics of price movements of on-line auctions of heterogeneous, hedonic products such as fine art (in our case, modern Indian art). It also provides a mechanism for us to analyze the relationship of factors that contribute to the heterogeneity on bid dynamics. This is particularly important for the auctioneers because it indicates the influence of various factors on bid values during an auction, which is vital for providing pre-auction estimates and selecting and optimizing lot sequence. The factors explored in this study are artist characteristics (established or emerging artist; prior sales history), art characteristics (size; painting medium—canvas or paper), auction design characteristics (opening bids; pre-auction estimates; position of the lot in the auction) and competitive characteristics (number of bidders; number of bids).

The rest of the paper is presented as follows. First, we discuss modern Indian art and describe the on-line auction data used in our analysis. Second, we explain the technique used to analyze the auction data of heterogeneous products and present various characteristics that may influence the dynamics of price formation. Third, we discuss the results obtained and the insights gleaned from them. Finally, we present the research challenges and opportunities that remain unexplored.

## 2. MODERN INDIAN ART AND ON-LINE AUCTION DATA USED

### 2.1 Modern Indian Art

Art auctions generated over $3 billion in sales in 2004 and the total world art market is estimated to be over $30 billion (www.artprice.com). The venerable auction houses, Christie's and Sotheby's, have dominated the art market since the 18th century. The advent of the Internet has introduced new players, who have adopted the new technology to threaten the dominance of the traditional players as well as to expand the market. One dramatic example is in the emerging market for modern Indian art. Although traditional auctions for modern Indian art have existed since 1995, it is only since 2000 that the market has exploded with values realized at auctions growing at a brisk 37.5% annually (coincidentally this is when SaffronArt.com, the source of our data, started its on-line auctions of modern Indian art). Modern Indian art in 2005 generated a total value in excess of $53 million with average lots selling for over $40,000. (This compares well with Latin American art, which is in a more mature stage in the auction market and generates about $50 million annually with average lot price of $50,000.) The top ten Indian artists sold 31% of the lots and contributed to 57% of the total value realized at auctions since 1995. Two of these artists are now ranked in the top 100 artists in the world based on their auction sales in 2005. A new set of emerging artists (the new trendsetters, typically born after 1955) have contributed 2% in value and 3% in lots, but are becoming increasingly popular, commanding increasingly higher prices. In the first half of 2005, SaffronArt.com, the on-line auction house of modern Indian art, sold almost as much Indian art as Christie's and sold more lots than any other auction house.

### 2.2 Data

We collected bidding data from SaffronArt.com, the leading on-line auction house of modern Indian art. The data for this study are the bid histories of 107 art lots (each lot typically is a unique piece of art, namely, a painting, a drawing or a sculpture; sometimes a lot can consist of a group of paintings or drawings; in our case, all lots have a single painting or a drawing) conducted in a three-day auction in December 2004. This auction uses an ascending-bid format with a fixed ending time and date set by the auction house, with the auction continuing until three minutes after the last bid arrival. Only bidders registered in advance for the specific auction may bid on the lots. The auction also uses a "proxy-bid" system similar to the one used by eBay, where the bids are automatically updated on behalf of the bidders.

In addition to the bid histories of the 107 lots, we collected information on artist characteristics (whether an emerging artist or an established artist; historical prices realized by the artist), art characteristics (size and medium of painting), auction design characteristics (pre-auction estimates; opening bid; order of the lot) and bidder competition characteristics (number of bidders; number of bids). The works of 48 artists were auctioned in this event, with an average of 2.4 lots per artist. The range of lots per artist varied between 1 and 13. The average realized price per lot was $19,406, and ranged from a low of $1,400 to a high of $132,222. The average price per square inch (this standardized measure of price is often an indicator of the value of the art and the artist) was $34.38 per square inch and ranged between $0.99 per square inch and $222.16 per square inch. There were a total of 127 unique bidders who participated in the on-line auction. The average number of bids per lot was 9.5, and it ranged between 2 and 23 bids. On average, there were 4 bidders participating in each lot, and ranged from 2 to 8 participants across the lots.

## 3. ANALYZING ON-LINE ART AUCTION DATA

### 3.1 Functional Data Analysis

The uniqueness of on-line auction data presents challenges that make the application of traditional econometric/regression methods difficult. The data consist of a sequence of bids placed over time. These bids do not arrive at evenly spaced intervals, precluding the use of traditional time series methods. Moreover, the bid dynamics change dramatically throughout the auction, particularly toward the beginning and the end of the auction. Functional data analysis (FDA) [25], which at its core is the analysis of curves rather than points, is well suited to analyze the bid dynamics of on-line auctions. Using this technique, we can recover the underlying price curves from which we can analyze the price dynamics (velocity and acceleration) in on-line auctions. Whereas traditional regression methods can be useful in modeling the final realized price in auctions, FDA provides the tools needed to model the price dynamics and determine the relationship of relevant strategic variables on price movement during the entire auction. This provides insights regarding how the momentum of the auction is being affected at various stages of the auction. To use FDA, we first smooth the bid data for each lot and recover the underlying price curves. Then we model the heterogeneity of these price paths using covariates of the art objects to provide insights on the relationship of these covariates to the price dynamics during the auction.

### 3.2 Recovering Price Curves

Whereas the bids arrive at irregular intervals and vary across lots in an auction, the raw bid data undergo some preprocessing steps. Let $t_{ij}$ be the time when the $i$th bid was placed, $i = 1, 2, \ldots, n_j$, for lot $j$ ($j = 1, \ldots, N$). Time $t_{ij}$ will vary for each lot due to the irregular arrival of bids. Although the auction event we are studying is a three-day auction, different lots have different closing times. Therefore, we standardize the auction time by scaling it between 0 and 1; thus $0 \leq t_{ij} \leq 1$. Now let $y_i^{(j)}$ be the bid amount placed at time $t_{ij}$. To realize the smooth price curve, we treat the $y_i^{(j)}$ values as noisy realization of the underlying continuous price function and reconstruct this function using appropriate smoothing techniques. [The smoothing here is done for a methodological reason, but there is an underlying conceptual reason why it makes sense to think of the price curve as continuous in the context of art auctions. In the case of a live auction, where the bid arrivals are usually compressed within 60 seconds (the time it typically takes to sell a lot), prices would appear to be continuous. Due to the nature of an on-line auction, bids arrive at different times over the duration of the auction, which is often spread over several days. Plotting the bids for lots where there is sufficient bidding activity over the duration of the on-line auction may appear choppy, but still would appear to be continuous.] We accommodate the irregular spacing of the bid arrivals by linearly interpolating the raw data and then sampling it at a common set of time points $t_i$, $0 \leq t_i \leq 100$, $i = 1, 2, \ldots, n$. Now the bids for the lots can be represented as vectors of equal length,

$$y^{(j)} = (y_1^{(j)}, \ldots, y_n^{(j)}),$$

where $y_i^{(j)} = y^{(j)}(t_i)$ represents the value of the interpolated bid sampled at time $t_i$. The bid values are then log transformed [Box–Cox normality test on bid price indicated $\lambda = 0$, suggesting a log transformation; a Box–Cox plot of the residuals before and after log transformation confirmed it as well] to reduce the high variability in bidding, especially near the end of the auction.

To recover the underlying price curves, we used penalized smoothing splines as suggested by Ramsay and Silverman [25] and Simonoff [28]. Smoothing splines provide small local variation and overall smoothness, and as such are a good choice for recovering the underlying curves. They also readily yield different derivatives of the target price curve, which is a desirable feature of the technique. Let us consider a polynomial spline of degree $p$,

$$f(t) = \beta_0 + \beta_1 t + \beta_2 t^2 + \beta_3 t^3 + \cdots + \beta_p t^p$$
$$+ \sum_{l=1}^{L} \beta_{pl} [(t - \tau_l)_+]^p,$$

where $\tau_1, \tau_2, \ldots, \tau_L$ is a set of $L$ knots and $u_+ = u I_{[u \geq 0]}$. The choice of $L$ and $p$ determines the departure of the fitted function from a straight line, with higher values resulting in a rougher $f$, which may result in a potentially better fit but a poorer recovery of the underlying trend. A roughness penalty function of the form

$$\text{PEN}_m = \int [D^m f(t)]^2 \, dt,$$

where $D^m f$, $m = 1, 2, 3, \ldots$, is the $m$th derivative of the function $f$, may be imposed to measure the degree of departure from the straight line. The goal is to find a function $f^{(j)}$ that minimizes the penalized residual sum of squares

$$\text{PENSS}_{\lambda,m}^{(j)} = \sum_{i=1}^{n} (y_i^{(j)} - f^{(j)}(t_i))^2 + \lambda \text{PEN}_m^j,$$

where the smoothing parameter $\lambda$ provides the trade-off between fit $[(y_i^{(j)} - f^{(j)}(t_i))^2]$ and variability of the function (roughness) as measured by $\text{PEN}_m$. [Sensitivity tests were performed with different values of $p$ (4, 5 and 6 were used) and $\lambda$ (14 different values between 0.001 and 100 were used). We found the model fit to be insensitive to different values of $p$ and $\lambda$. However, the root mean squared error for the model was the lowest with $p = 4$ and $\lambda = 0.1$. Thus, we use these smoothing parameters to recover the price curves.] We used the monospline module developed by Ramsay [24] to minimize $\text{PENSS}_{\lambda,m}^{(j)}$.

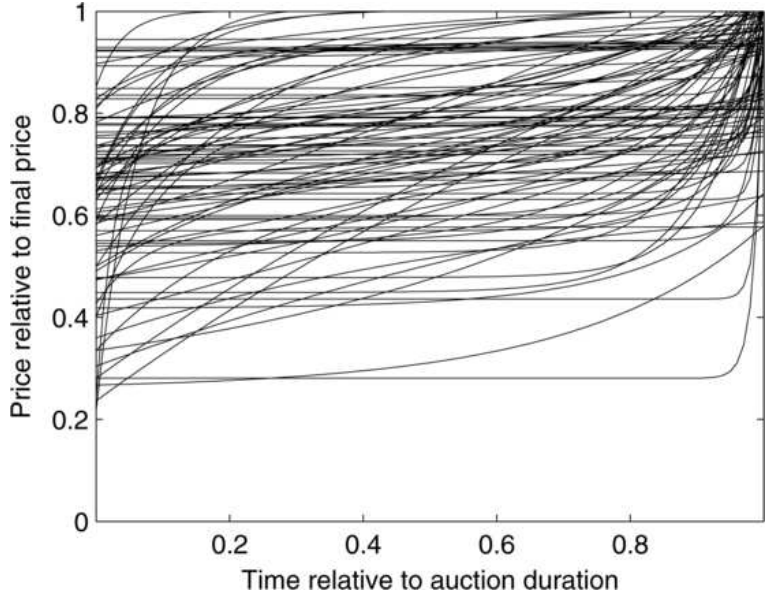

Fig. 1. *Smoothed price curves of* 107 *auction lots.*

Let $f_j(t)$ be the smoothing spline for lot $j$. We estimate a smoothing spline for each auction lot. Such smoothing splines yield several derivatives of the bid path curve that provide us with a detailed look at the underlying dynamics, such as velocity [first derivative $f'_j(t)$] and acceleration [second derivative $f''_j(t)$]. Figure 1 shows the plots of the smoothed (and log transformed) curves of the 107 lots. Since different lots have different realized prices, we scaled the price curves between 0 and 1, where each data point represents the value of the item as a fraction of the final realized value of the item at that instant in the auction. These will be our price curves, which will be used in further analysis in place of the original data. Figure 2 shows the different price paths of some select auction lots to illustrate the variety of price curves that is evident in our on-line auction. The price path of lot 81, a painting by an emerging artist, N. Sharma, shows that the price starts at about 30% of its eventual realized price ($12,094) and moves in a linear fashion as the auction progresses (see Table 1). Lot 1, a painting by J. S. Ali, shows an increasing curve with an increasing rate. It starts as a flat price path until the last quarter of the auction, when it moves rapidly toward its realized value ($7,794). The price path of lot 16, a painting by one of the leading established artists, R. Kumar, also shows a curve at an increasing rate, but the price starts out at about 80% of the realized value and shows a steady path until the end of the auction. Lot 10, a work on paper by J. Sabavala, shows an increasing curve at a decreasing rate, where most of the price movement appears to be during the first quarter of the auction.

### 3.3 Functional Regression

As we see from the previous discussion, the estimated price curves showed a great degree of variability. It would be of interest to investigate the factors that explain these patterns. Unlike standard regression models, where predictor and explanatory variables are scalars or vectors, functional regression allows one to have these variables take on a functional form. For example, the response variables in our case are the price curve $f(t)$, $f'(t)$ (velocity) and $f''(t)$ (acceleration) that capture the price formation process during the auction. Potential explanatory variables are characteristics of the lot, such as the type of painting, size of the painting, artist reputation, opening bid, and will be described in great detail below. To capture the effects of these variables on the price dynamics, we run a regression for each time period (1–100) for data from all the lots ($n = 107$). The parameter estimates associated with each explanatory variable are then plotted along with confidence bands to indicate the effect and its significance over the entire auction.

TABLE 1
*Characteristics of the lots presented in Figure 2*

| Lot | No. of bids | No. of bidders | Name | Pre-auc. low est. | Pre-auc. high est. | Realized value | Length of the art work (inches) | Width of the art work (inches) | Medium | Artist | Artist type |
|---|---|---|---|---|---|---|---|---|---|---|---|
| 1 | 12 | 4 | Krodhit | $4,000 | $5,000 | $7,794 | 40.5 | 68.5 | Oil on canvas | Ali, J. S. | Others |
| 10 | 5 | 4 | Untitled | $6,340 | $8,560 | $7,794 | 17 | 14 | Charcoal on paper | Sabavala, J. | Others |
| 16 | 11 | 5 | Untitled | $26,670 | $31,120 | $46,225 | 36 | 36 | Oil on canvas | Kumar, R. | Established |
| 81 | 15 | 5 | Fat, F*** and Forty | $5,000 | $6,000 | $12,094 | 30 | 71 | Acrylic on canvas pasted on board | Sharma, N. | Emerging |

The following sets of explanatory variables are considered in explaining the price curves of lots in the on-line auction of modern Indian art.

*Auction house/seller design characteristics.* The seller, in this case, the on-line auction house, has the choice of several design issues, such as the opening bid level (set usually as a percentage of the pre-auction low estimate), the value of a hidden reserve price, the position of a lot in the auction and the duration of the auction. Whereas all the lots in this auction had a hidden reserve and all the lots ended at the end of the three days, these two factors will not be considered in our analysis. Unlike eBay auctions, pre-auction low and high estimates are provided for each lot along with a suggested opening bid. Evidence of the relationship of the opening bid and the pre-auction estimates on prices has been found by some researchers. Czujack and Martins [12], using traditional auction data, found that the presale estimates are a good predictor of the realized price. Bajari and Hortaçsu [3] found that low opening bids attract more bidders. Bapna, Jank and Shmueli [8] found that higher opening bids result in lower auction momentum and slower price formation. Furthermore, the importance of opening bid should decline as the auction progresses, since other characteristics of the artwork, like the quality, size and genre, may become more relevant in influencing the price formation process. Therefore, we anticipate that price level will be positively associated with opening bid and that this association will diminish as the auction progresses.

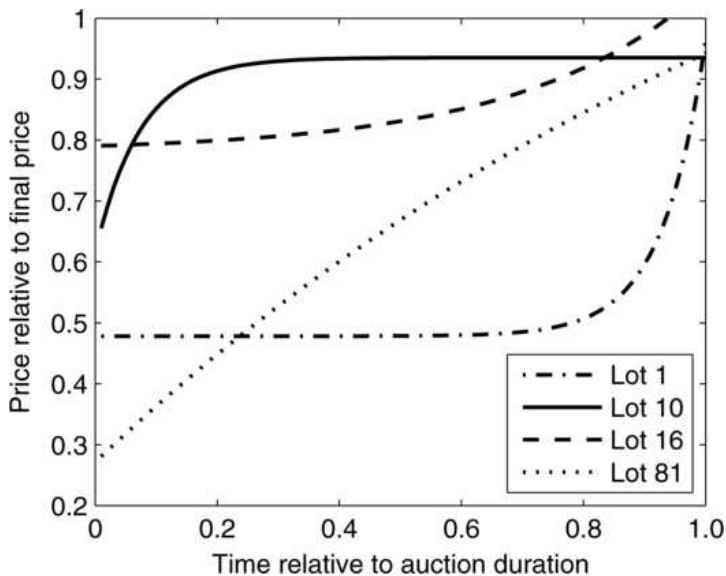

FIG. 2. *Auction lots with different price paths.*

In a typical auction, about 100 lots are sold and each one is assigned a sequential lot number that indicates the order in which the lot will be sold. In the current auction, lots are broken up into groups of 20–25 each. There were five such groups. The auction ending time is staggered for these groups. For example, the auction ends for the first group of lots at 9:30 am, the ending time for the second group of lots is 30 minutes later, and so on. Although the arrangement of lots into groups and the ordering of these groups is under the control of the auction house, there does not appear to be a system or rationale by which these decisions are made. For instance, lots of the same artists are often placed together or lots with high estimates are placed at the beginning of the auction; sometimes young or emerging artists are placed at the end of the auction. Beggs and Graddy [10], talking to representatives at Christie's Impressionist and Modern Art department, reported that lots are "loosely ordered by date," with earlier works being placed earlier in the auction. Examining the way that the lots and groups are ordered in this auction shows that lots of established artists are placed in the group that ended first and the emerging artist lots were placed in the group that ended last, with no evident pattern of placement or ordering of the groups ending in the middle. Very little empirical evidence exists as to the effect of lot order on prices. Evidence from research on sequential auctions (auctions of similar items in the same auction) is quite ambiguous: some researchers found a price decline for similar lots over an auction and others found a price increase [21]. Nunes and Boatwright [22] found some order effects in automobile auctions. They found that a high bid or price paid for a previous automobile significantly (and positively) affects the price of the focal car. Beggs and Graddy [10] report that an auction house might place a watercolor immediately after the sale of a major oil painting by the same artist to benefit from the "excitement" created by the sale of the oil. Although these results indicate that order may play a role, little evidence exists as to its effects. Based on the way the lots are ordered in this auction, we anticipate that the price formation would be affected in a negative way, with lots that end earlier showing higher current prices.

*Competition characteristics.* The level of competition in an auction may be captured by the number of bidders. Competition among bidders plays an integral role in price formation. As the number of bidders increases, the competition among them also increases, resulting in a higher price. Bapna, Jank and Shmueli [8] found the number of bidders to be positively associated with the current price of the item. Furthermore, it is observed that, typically, the middle of the auction experiences a smaller amount of bidder participation as compared to the early and later stages of the auction. Bidders generally utilize this time to scrutinize the auctioned item or just simply wait to see how other bidders behave. Therefore, it would be interesting to see how this competition characteristic affects the on-line art auction's price formation. We anticipate that the number of bidders has a significant positive relationship with price levels. (Another interesting competition characteristic is the current number of bids for an item, which also signals the level of competition among the bidders. An increase in the number of bids signals higher competition among the bidders and thus may result in a higher price. In our study, number of bids is highly correlated with number of bidders and hence adds little additional explanatory power.)

*Artist characteristics.* Artist characteristics like artist reputation and previous auction history tend to provide a basis for the bidders to estimate the value of the item. For example, established artists are more highly reputed and recognized for their work than are emerging artists. Most of their works have been resold many times in the market and, thus, the bidders know their values. Therefore, the works of established artists present a low risk investment opportunity for the bidders, but are bound to be more expensive. On the other hand, emerging artists are new to the art market, their works have not reached the saturation point, are more affordable and therefore are likely to attract more bidder and price activity, and hence possess greater energy during the auction. We have classified the artists represented in this auction into three categories. Established artists are those born in the first quarter of the 20th century who started the modernist art movement in India. These are artists like M. F. Husain and F. N. Souza, who are currently ranked in the top 100 artists of the world based on their auction sales. The second group is the emerging artists, who are artists born after 1955, are setting new trends in the Indian art scene and are becoming increasingly popular. Examples of emerging artists are Atul Dodiya, Shibu Natesan, Baiju Parthan and Chitravonu Mazumdar. All artists who do not fall into these two categories are classified as "other" artists. Our data set contains works of established artists ($n = 33$), emerging artists ($n = 20$) and other artists ($n = 54$) who fall in neither of these categories. In our analysis, we study the relationship between artist type and price formation. We anticipate that established artists have a positive effect and emerging artists have a less positive effect on price formation during the auction.

Historical market information, like the number of lots sold by the artist in the previous year and the average value realized (price per square inch) for the artist's works in the previous year, provides a signal to the bidders about the market value of the artists. Artnet.com and Artprice.com are two of several sources that the bidders may consult prior to and during the auction with regard to past auction activity of the artists in whom they are interested. Sotheby's has recently teamed up with Artnet.com and started to provide detailed historical information about the artists on its web site before an auction. Bidders typically use this information to choose the lots to bid on and the amount to bid. It will be interesting to evaluate how these pieces of information affect the price formation process in on-line art auctions. We expect them to have a significant positive effect on the price at least during the beginning of the auction, with greater past activity and higher value realized for an artist providing a positive signal to bid a higher price. (Due to the high correlation between the number of lots sold and the prior value realized for the artists, we will use only historical prices realized in the model.)

TABLE 2
*Summary statistics of the explanatory variables*

|  |  | Mean (SD) | Median | Min. | Max. |
|---|---|---|---|---|---|
| Artist characteristics | Average price/sq. inch realized in previous year's auctions | $36 (42) | $10.29 | $1 | $222 |
|  | Lots sold by the artist in previous year's auctions | 19.33 (19) | 11 | 1 | 68 |
| Auction design characteristics | Pre-auction low estimate in $ | $13,018 (16,131) | $6,230 | $780 | $100,000 |
|  | Opening bid in $ | $11,195 (13,943) | $4,975 | $622 | $87,500 |
| Competition characteristics | No. of unique bidders/lot | 4.06 (1.64) | 4 | 2 | 8 |
|  | No. of bids/lot | 9.504 (5.159) | 8 | 2 | 23 |
| Art characteristics | Size of the art work in sq. inches | 1,273 (1.304) | 822 | 52.25 | 5,796 |

*Art characteristics.* The functional attributes of fine art may play a role in price formation during the auction. Larger paintings tend to command higher absolute prices than smaller ones, but may be a better value based on price per square inch. In traditional auctions, Czujack, Flores and Ginsburgh [11] found that size/area of the painting has a negative impact on price indices of art over time. Beggs and Graddy [10], in their examination of auction prices of Impressionist and Modern Art, found that, in general, the size of the paintings in terms of length and width is positively associated with bid prices, but the width of the painting consistently had a larger effect on bid prices than length. Also, paintings on canvas command a higher price than works on paper. However, we are not aware of any study that has examined the effect of size and art medium on the price formation process in an auction.

We defined the following characteristics and used them in the functional regression of log(price) (except for the dummy variables, all variables have been log-transformed based on the Box–Cox normality tests):

$x_{j1} = $ log(price per square inch of lots sold by the artist of lot $j$ in previous year's auctions);

$x_{j2} = $ dummy variable to indicate if lot $j$ belonged to an established artist ( $= 1$ if an established artist);

$x_{j3} = $ dummy variable to indicate if lot $j$ belonged to an emerging artist ( $= 1$ if an emerging artist);

$x_{j4} = $ log(opening bid of lot $j$);

$x_{j5} = $ log(position (#) of lot $j$ in the auction) ( $= 1$ if the group in which it is placed ended first; $= 2$ if the group in which it is placed ended second, etc.);

$x_{j6} = $ log(size (area in square inches)) of lot $j$;

$x_{j7} = $ dummy variable to indicate medium of lot $j$ ( $= 1$ if canvas; $= 0$ if paper);

$x_{j8} = $ log(number of bidders for lot $j$ at time $t$).

TABLE 3
*Correlation matrix of explanatory variables*

|  | Previous price per sq. inch of the artist ($x_1$) | Established artist ($x_2$) | Emerging artist ($x_3$) | Opening bid ($x_4$) | Lot position ($x_5$) | Size of the painting ($x_6$) | Canvas ($x_7$) | No. of bidders ($x_8$) |
|---|---|---|---|---|---|---|---|---|
| Previous price per sq. inch of the artist ($x_1$) | 1 | | | | | | | |
| Established artist ($x_2$) | 0.55** | 1 | | | | | | |
| Emerging artist ($x_3$) | −0.41** | −0.38** | 1 | | | | | |
| Opening bid ($x_4$) | 0.59** | 0.38** | −0.23** | 1 | | | | |
| Lot position ($x_5$) | −0.50** | −0.64** | 0.69** | −0.37** | 1 | | | |
| Size of the painting ($x_6$) | −0.27** | −0.25** | 0.50** | 0.17 | 0.33** | 1 | | |
| Canvas ($x_7$) | 0.22* | −0.04 | 0.10 | 0.41** | 0.00 | 0.46** | 1 | |
| No. of bidders ($x_8$) | −0.19* | 0.00 | 0.18 | −0.05 | 0.06 | 0.22** | 0.09 | 1 |

*Correlation is significant at the 0.05 level (two-tailed).
**Correlation is significant at the 0.01 level (two-tailed).

Dummy variables are used to identify the medium and the artists. The variables that deal with number of bidders and number of bids are created differently to capture their changing nature over the auction. The dynamic covariates are created by using an evenly spaced grid over the duration of the auction, $t_1, t_2, t_3, \ldots, t_n$ (in this case $n = 100$), and the number of bids and number of bidders are assigned for each period $t_i$. (Although we created the dynamic measure for number of bids, we did not use this variable in our final analysis because it was highly correlated with number of bidders.) Summary statistics of the explanatory variables and their correlations are presented in Table 2 and Table 3, respectively.

## 4. RESULTS

The results of the functional regression are displayed in the form of estimated parameter curves (solid lines) in Figure 3. The summary of findings is presented in Table 4. A plot of the confidence bands (gray lines) around the curve is also displayed to allow us to draw inference (at 95% confidence level) about the significance of the estimated curve. The results from the functional regression with the covariates are presented next. We will present and discuss the parameter curves of price $f(t)$ and velocity $f'(t)$. The regressions of price acceleration curves $f''(t)$ were noisy and yielded no useful insights, and hence are not presented here. (To investigate the influence of outliers on our results, we identified observations that are more than 2.5 standard deviations away on any one of the following variables: opening bid, size of the lot and price per square inch of the lots sold in the previous year for the artist. This resulted in a sample size of 100. The results obtained from this sample were not different from those of the full sample. A suggestion made by a referee prompted us to analyze the model by excluding lots with only two bidders or bids to see if this influenced the results. Whereas the results were not different, we chose to report the results with all observations.)

TABLE 4
*Summary of findings*

|  | **Price level** | **Price dynamics (velocity)** |
|---|---|---|
| | *Auction design characteristics* | |
| Opening bid | Opening bid is positively related with price level. This relationship declines as the auction progresses. | Higher opening bids associated with steady price increases during most of the auction, with the rate of increase diminishing at the end of the auction. |
| Lot position | Lot position negatively related to price level, indicating later lots have lower price levels. This effect is negative throughout the auction. | Lot position is associated with slower price increases throughout the auction. |
| | *Competition characteristics* | |
| Number of bidders | Not significant, but the directionality of the relationship suggests positive influence of the number of bidders on price levels at the beginning of the auction, showing a decline toward the end of the auction. | Not significant, but the effect suggests a slower increase at the beginning of the auction and a faster price increase during the middle of the auction. |
| | *Artist characteristics* | |
| Established artists | Established artists show higher price levels at the beginning of the auction. This effect declines steadily toward the end of the auction. | Faster price increases as the auction progresses, picking up momentum toward the end of the auction. |
| Emerging artists | Emerging artists show lower price levels at the beginning of the auction. This effect diminishes steadily toward the end of the auction. | Slower price increases throughout the auction. |
| Historical value of artists (price/sq. inch) | Past value of the artist affects price level positively at the beginning of the auction. This relationship is still strong, but declines toward the end of the auction. | Faster price increases throughout the auction, getting slower in the middle and the end of the auction. |
| | *Art characteristics* | |
| Size of the painting | Price levels are lower for larger paintings at the beginning of the auction. | Slower price increases at the beginning of the auction for larger paintings. |
| Canvas | Not significant. | Not significant. Faster price increases at the end of the auction for canvas paintings. |

*Auction house/seller design characteristics.* The estimated parameter curve for opening bid is positive throughout the auction, implying a positive relationship between the size of the opening bid and the current price. The confidence bounds (at 95% confidence levels) indicated by the gray lines signify that the relationship is statistically significant throughout the auction except at the very end of the auction (see Figure 3). The parameter declines over the period of the auction, implying that the effect of the size of the opening bid on price formation declines over the auction, with the greatest effect during the beginning of the auction. So the opening bid provides the most information during the start of the auction but fails to be useful in explaining the price at the end of the auction. The effect of the size of the opening bid on the price dynamics (velocity) during the auction was not significant at the beginning and at the end of the auction, but was positive and significant during the middle of the auction. Contrary to what Bapna, Jank and Shmueli [8] found in eBay auctions—that high opening bids provide less energy resulting in slower *price increases* (rate of change in price) during the auction— the effect of opening bid price on velocity in the current setting appears to be positive, indicating a steady price increase throughout the auction. Unlike at eBay, where the sellers have the option to set the opening bid as low as $0.01, the opening bid in the current auction is pegged to the pre-auction low estimate and provides much less variation in this variable. (Opening bid by this on-line auction house is set as a percentage of the pre-auction low estimate. Here it was about 85% of the low estimate and ranges from 73 to 96%.) A closer examination of the data shows that lots with opening bids of less than $5,000 had more bidders and almost twice as many bids as lots with higher opening bids. The findings from the eBay auctions may provide some strategic guidance in auction design in terms of the choice of opening bid. Opening bids are used to protect the seller. However, whereas most lots have a reserve below which the seller will not sell, low opening bids could be used as a strategic option by the seller to provide energy during early periods of the on-line auction. People familiar with live auctions will note that the auctioneer typically starts the bids at about 50% of the low estimate, but on many occasions brings down the opening bid (way below to about 10% of the low estimate) to encourage bidding. This invites bidders to bid and the lot may often be sold for substantially higher than the low estimate.

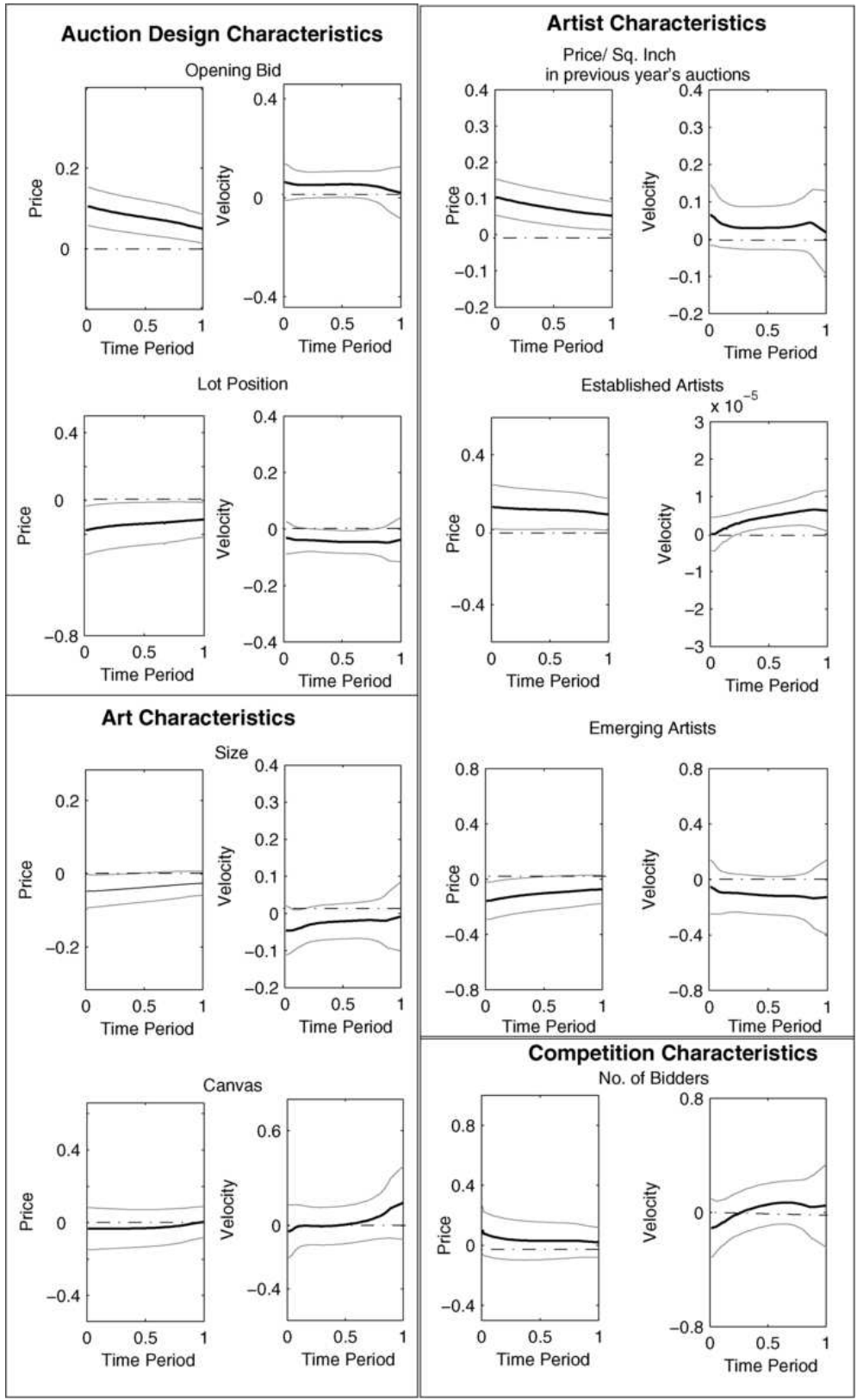

Fig. 3. *Estimated parameter curves for auction design characteristics, artist characteristics, art characteristics and competition characteristics.*

The parameter curve for the position/order of the lot in the auction suggests a significant negative relationship between the lot order and the current price. Early lot placement yields a higher current price. This effect diminishes as the auction progresses. The parameter curve associated with price velocity, although not significant at the beginning and at the end of the auction, is negative, indicating a steady but slower price increase during the middle of the auction. The average pre-auction low estimate of the lots in group 1 (which ended first) was $21,000, and for lots in group 5, which ended last, was $5,000. From this, we observe that the auction house's design decision on ordering the lots is partly based on value. This is a simple and easily implementable design. Will alternative ordering enhance the dynamics of the auction? The findings of Beggs and Graddy [10] in art auctions and Nunes and Boatwright [22] in automobile auctions suggest some interesting options to explore. Will the placement of an important artist's works spread over the auction enhance the prices? Will an important painting by an important artist followed by a painting by an emerging artist provide greater energy and yield higher prices? Examining such strategies using FDA will provide valuable insights on the price formation and the dynamics of price during an auction. Such insights are valuable because they can help auctioneers order the lots to maximize prices realized.

*Competition characteristics.* Although the current number of bidders was not significant in this study, the directionality was similar to that found by Bapna, Jank and Shmueli [8]. Small sample size and low variability of these two variables in our data might have contributed to the insignificant result here. With increasing attention paid to modern Indian art, bidder activity is only going to increase in the future, which will strengthen the effect of competition observed here. [A recently concluded on-line auction of modern Indian art (December 2005) showed that the total number of bidders has increased to 299, a 135% increase since December 2004.] The effect was positive throughout the auction, with its effect diminishing toward the end of the auction. The parameter curve associated with price velocity indicated slower price increase at the beginning of the auction, but the rate of change steadily increased during the middle of the auction.

One other way to enhance competition is to provide access to information on bidders and bidder activity during the auction. eBay provides information on all bidders and their activity in an easily accessible manner to all users. The on-line auction house studied here provides bidder information on a lot by lot basis with regard to who the bidders are (anonymous bidder identifications are assigned before the start of the auction), what the bid amounts were and when the bids were placed. However, it does not provide information on what different lots the bidders are actively bidding on and how much they are bidding on other lots in an easily accessible manner. The provision of such information may enhance competition and influence the price dynamics positively during the auction.

*Artist characteristics.* The artist reputation and previous auction history is captured by several measures that provide some interesting insights. The relationship between established artists and price formation is positive and significant, and has the strongest influence at the beginning of the auction. The prices of these artists tend be fairly high early in the auction (recall the price path of lot 16 of R. Kumar, an established artist presented in Figure [2](#)). This effect shows a steady decline toward the end of the auction. An opposite and symmetric effect is seen for emerging artists. The effect is negative and significant at the beginning of the auction, implying that prices realized at the beginning of the auction are lower than those for established artists. The effect becomes insignificant at the end of the auction, suggesting the prices realized for emerging artists tend to catch up a bit with those obtained for established artists. The parameter curve associated with price velocity for established artists is positive for most of the auction (except at the beginning and at the end of the auction), indicating faster price increases during the middle of the auction. This momentum slows down by the end of the auction.

The relationship between the auction prices realized by the artists' works in the previous year and current price level is positive and significant throughout the auction, implying that the artists with established reputation and record attract higher bids. In both of these cases, however, the strength of this relationship decreases steadily as the auction progresses. This declining effect implies that the impact of the artist's reputation on price formation decreases throughout the auction, and information contained in the reputation measures has greater signaling value in the early stages of the price formation process and becomes less useful as the auction progresses. As the auction progresses, bidders potentially receive other relevant information from other bidders and possibly from other sources, making the artist reputation information not as critical in their bid process. The parameter curve associated with price velocity was positive but not significant during most of the auction. The positive effect suggests that the signal provided by artists' past records results in faster price increases during the auction but may be more pronounced at the beginning of the auction. These results suggest an opportunity for the auction house to provide comprehensive historical auction information about the artists to all potential bidders, thus leveling the playing field. Sotheby's auction house, as we reported earlier, is doing exactly this.

*Art characteristics.* The medium of painting (canvas) did not have a significant effect on price formation during the auction. However, the size of the paintings had a significant negative effect on current price at the beginning of the auction, indicating that current price as a percentage of final realized value will be lower for larger paintings at the beginning of the auction. Examining the parameter curve for the velocity indicates a slower price increase for larger paintings. When we regressed the different dimensions of the painting (length and width) on the price curves (instead of size), we found that the width of the painting had a negative effect on current price but the length did not. This is contrary to the findings reported by Beggs and Graddy [10] in their examination of auction prices of Impressionist and Modern art. They found that, in general, both length and width were positively associated with bid prices, but the width of the painting consistently had a larger effect on bid prices than length. The reason for the difference could be due to the nature of the art examined here. Modern Indian art is still an emerging art form, with a large number of lots coming from young emerging artists compared to the established nature of Modern and Impressionist art.

## 5. CONCLUSIONS AND FUTURE DIRECTIONS

Standard statistical procedures are adequate to analyze and model price data that are obtained from traditional auctions. However, these methods are not adequate to deal with the dynamics of the price formation data and the bidder data that are available through on-line auctions. In this paper, we illustrate how functional data analysis can be useful to understand the price formation of on-line fine art auctions. FDA provides a way to see the effect of relevant strategic variables on price movement during the entire auction rather than simply at the end price. This provides insights regarding how the momentum of the auction is being affected at various stages of the auction and provides options for the on-line auctioneer in the design of the auction (ordering the lots, determining the opening bid) to maximize realized prices.

Our analysis suggests that the opening bid positively affects on-line auction price levels of art at the beginning of the auction and its effect declines toward the end of the auction. The order in which the lots appear in an art auction has a negative effect on current price level, with this relationship decreasing toward the end of the auction. This implies that lots that appear earlier have higher current prices during the early part of the auction, but that effect diminishes by the end of the auction. Established artists show a positive relationship with price level at the beginning of the auction. Reputation or popularity of the artists and their investment potential as assessed by previous history of sales have a positive effect on prices at the beginning of the auction. The medium (canvas or paper) of the painting did not affect art auction price levels, but size of the painting was negatively related to the current price during the early part of the auction.

Important implications for on-line art auction houses and bidders can be drawn from this study. It is clear from the findings that effects on price during an auction do vary over the auction, with some of these effects being strong at the beginning of the auction (such as the opening bid and historical prices realized). In some cases, the rate of change in prices (velocity) seems to increase at the end of the auction (for canvas paintings and paintings by established artists). The potential use of opening bid as a strategic variable to induce more energy in the auction should be seriously considered. The order in which lots are placed is another important strategic variable that the on-line auction house can use to maximize realized prices. Providing easy access to information on artists' historical auction records and on real time bidder activity will spur competition and help bidders make more informed choices during the auction.

The bases of the current research are the smoothed price curves of 107 lots in an on-line art auction. This provided us with some interesting insights on the magnitude of current price, rate of change of these prices and the effects of covariates during an auction. In the future, research should shed light on the other dynamics of price change, namely acceleration of price. Also, the results obtained here can be generalized after replicating the findings from data obtained from other on-line art auctions.

The data in on-line art auctions are inherently hierarchical, with the lots in the auction (107) nested within artists (48). So the use of a mixed model, which allows for the dependence among lots belonging to the same artist, may be appropriate. In the future, developing such a mixed model and comparing its results with those obtained from FDA will provide a better understanding of on-line art auctions.

The availability of dynamic auction data and the appropriate statistical tools to analyze them provide some interesting opportunities and challenges to researchers in this area. Following the on-line art auctions of a particular genre (e.g., modern Indian art) over time and collecting the auction dynamic data provides us with some fruitful areas for future research.

Do auction dynamics change as markets or products evolve? Modern Indian art is seeing an explosive growth at this time. What kind of auction dynamics can we expect and how will the various factors examined here affect these dynamics in more mature markets or when the modern Indian art market reaches its mature stage? As we follow this and other markets over the next few years, answers to these important questions may be found.

Several pieces of work of some artists keep appearing in on-line auctions over time. For example, almost a dozen works of the well known Indian artist M. F. Husain are routinely sold at each auction. It would be interesting to see how the auction dynamics for specific artists evolve over time. In addition to the changes in the final prices realized for the artist's work, one could analyze whether the underlying price formation process is changing as well. As the artist's value increases over time, does the price formation dynamic follow differing patterns? Another important area of inquiry is to examine the auction dynamics of paintings that appear at auctions repeatedly over time. Mei and Moses [19] evaluated the investment returns of art, based on an analysis of several thousand paintings which came up for sale over a 100-year period. As we capture data from more on-line art auctions, it will be interesting to analyze the price dynamics of a given painting that comes to auction repeatedly over time.

There is some evidence in the economics literature on art auctions of a declining price anomaly, where the prices of similar lots tend to decline over a short period of time. For example, the same Picasso print that is sold at a second auction within a short period after the first realizes a lower price. This phenomenon, which was first shown by Ashenfelter [1] and subsequently by other researchers, is commonly attributed to market inefficiencies. Very little work has been done in the on-line auction arena regarding this subject. Developing FDA models that incorporate data from a panel of auctions will help in understanding this anomaly fully.

Analysis in this research suggests interesting insights about the lot placement/order in an auction. This design element has not been studied in any systematic way. Nunes and Boatwright's [22] finding of a positive effect of a lot which sold for a high price on the lot immediately following it (in the case of automobile auctions) is fascinating. Future FDA models that incorporate lead-in and lead-out lot characteristics may shed light on this phenomenon and provide important design insights for auction houses in the placement of lots in an auction.

The community of bidders for fine art of a specific genre is usually small (about 200–250 for modern Indian art according to current estimates), but potentially heterogeneous in their bidding behavior. Bapna, Goes, Gupta and Jin [6] find significant but stable heterogeneity among bidders in on-line auction markets. See Dass and Reddy [13] for a classification of the bidder types in on-line art auctions. Unlike in eBay auctions, bidders in on-line art auctions need not have the same identity, which makes it difficult to follow their purchases over time, but if this information is made available by the auction house, characteristics of the bidders (number of lots; value and type of art bought before) can be incorporated, which will enrich the FDA models further.

Ku, Malhotra and Murnighan [18] found what they called competitive arousal or auction fever in several auctions that they studied, where bidders appear to bid substantially higher than their reservation price. On-line art auctions may also provide a fertile ground for understanding the dynamic learning behavior of bidders as the auction progresses. Do bidders change their bidding behavior during the auction? Do they learn from the lots that they won or lost and alter their bidding behavior on subsequent lots in the same auction? Creative models to understand such behavior should be pursued in the future.

Although most auctions tend to have some elements of private and common value elements, the setting we chose to study here, namely on-line art auctions, could be considered predominantly a private value setting. The results we discussed here have implications for such auctions. How will these implications transfer to a common value setting? What conclusions could we draw with respect to bidder and seller surplus? These important questions need to be addressed in future research.

## ACKNOWLEDGMENTS

We would like to thank Professor Wolfgang Jank and Professor Galit Shmueli of the Robert H. Smith School of Business at the University of Maryland for advice and guidance throughout this research. We would also like to thank Dr. Umesh Gaur and Professor Vanessa Patrick for comments and suggestions. Research support from the Coca-Cola Center for Marketing Studies is acknowledged.

## REFERENCES

[1] ASHENFELTER, O. (1989). How auctions work for wine and art. *J. Economic Perspectives* **3**(3) 23–36.
[2] ASHENFELTER, O. and GRADDY, K. (2002). Art auctions: A survey of empirical studies. CEPR Discussion Paper 3387. Available at www.econ.ox.ac.uk/members/kathryn.graddy/research/auctionsurvey.pdf.
[3] BAJARI, P. and HORTAÇSU, A. (2004). Economic insights from Internet auctions. *J. Economic Literature* **42** 457–486.
[4] BAPNA, R., GOES, P. and GUPTA, A. (2001). Comparative analysis of multi-item online auctions: Evidence from the laboratory. *Decision Support Systems* **32** 135–153.
[5] BAPNA, R., GOES, P. and GUPTA, A. (2003). Analysis and design of business-to-consumer online auctions. *Management Sci.* **49** 85–101.
[6] BAPNA, R., GOES, P., GUPTA, A. and JIN, Y. W. (2004). User heterogeneity and its impact on electronic auction market design: An empirical exploration. *MIS Quarterly* **28** 21–43.
[7] BAPNA, R., GOES, P., GUPTA, A. and KARUGA, G. (2002). Optimal design of the online auction channel: Analytical, empirical and computational insights. *Decision Sciences* **33** 557–577.
[8] BAPNA, R., JANK, W. and SHMUELI, G. (2004). Price formation and its dynamics in on-line auctions. Working paper RHS-06-003, Smith School of Business, Univ. Maryland. Available at ssrn.com/abstract=902887.

[9] BAUMOL, W. J. (1986). Unnatural value: Or art investment as floating crap game. *Amer. Economic Review* **76**(2) 10–14.
[10] BEGGS, A. and GRADDY, K. (1997). Declining values and the afternoon effect: Evidence from art auctions. *RAND J. Economics* **28** 544–565.
[11] CZUJACK, C., FLORES, JR., R. and GINSBURGH, V. (1996). On long-run price comovements between paintings and prints. In *Economics of the Arts: Selected Essays* (V. Ginsburgh and P. Menger, eds.) 85–112. North-Holland, Amsterdam.
[12] CZUJACK, C. and MARTINS, M. F. O. (2004). Do art specialists form unbiased pre-sale estimates? An application for Picasso paintings. *Appl. Economic Letters* **11** 245–249.
[13] DASS, M. and REDDY, S. K. (2005). Bidder heterogeneity in online art auctions. Working paper, Terry College of Business, Univ. Georgia.
[14] ELMAGHRABY, W. (2004). Auctions and pricing in E-Marketplaces. In *Handbook of Quantitative Supply Chain Analysis: Modeling in the E-Business Era* (D. Simchi-Levi, S. D. Wu and Z. M. Shen, eds.) 213–246. Kluwer, Norwell, MA.
[15] JANK, W. and SHMUELI, G. (2006). Functional data analysis in electronic commerce research. *Statist. Sci.* **21** 155–166.
[16] JAP, S. (2002). Online reverse auctions: Issues, themes, and prospects for the future. *J. Academy of Marketing Science* **30** 506–525.
[17] JAP, S. (2003). An exploratory study of the introduction of online reverse auctions. *J. Marketing* **67**(3) 96–107.
[18] KU, G., MALHOTRA, D. and MURNIGHAN, J. K. (2005). Towards a competitive arousal model of decision-making: A study of auction fever in live and Internet auctions. *Organizational Behavior and Human Decision Processes* **96** 89–103.
[19] MEI, J. and MOSES, M. (2002). Art as an investment and the underperformance of masterpieces. *Amer. Economic Review* **92** 1656–1668.
[20] MITHAS, S., JONES, J. L. and MITCHELL, W. (2005). When do firms use Internet-enabled reverse auctions? The role of asset specificity, product specialization and non-contractibility. Working paper, Smith School of Business, Univ. Maryland.
[21] NATZKOFF, P. (2001). Sequential auctions: Survey of the literature and empirical analysis of auctions of raw wool in Australia. M. Phil. thesis, Univ. Oxford.
[22] NUNES, J. C. and BOATWRIGHT, P. (2004). Incidental prices and their effect on willingness to pay. *J. Marketing Research* **41** 457–466.
[23] PESANDO, J. E. (1993). Art as an investment: The market for modern prints. *Amer. Economic Review* **83** 1075–1089.
[24] RAMSAY, J. O. (2003). Matlab, R, and S-PLUS functions for functional data analysis. Online resource. Available at ftp://ego.psych.mcgill.ca/pub/ramsay/FDAfuns.
[25] RAMSAY, J. O. and SILVERMAN, B. W. (1997). *Functional Data Analysis*. Springer, New York.
[26] SHMUELI, G. and JANK, W. (2006). Modeling the dynamics of online auctions: A modern statistical approach. In *Economics, Information Systems and E-commerce Research II: Advanced Empirical Methods* **1** (R. Kauffman and P. Tallon, eds.). Sharpe, Armonk, NY. To appear.
[27] SHMUELI, G. and JANK, W. (2005). Visualizing online auctions. *J. Comput. Graph. Statist.* **14** 299–319. MR2160815
[28] SIMONOFF, J. S. (1996). *Smoothing Methods in Statistics*. Springer, New York. MR1391963
[29] SPANN, M. and TELLIS, G. J. (2006). Does the internet promote better consumer decisions? The case of name-your-own-price auctions. *J. Marketing* **70**(1) 65–78.